\documentclass[twoside,11pt,leqno]{amsart}
\usepackage{amssymb}
\swapnumbers
\theoremstyle{plain}
\newtheorem{thm}{Theorem}[section]
\newtheorem{prop}[thm]{Proposition}
\newtheorem{cor}[thm]{Corollary}
\theoremstyle{definition}
\newtheorem{defn}[thm]{Definition}
\newcommand{\proofof}[1]{\end{#1}\begin{proof}}
\newcommand{\emphdef}{\textit}
\numberwithin{equation}{section}

\newcommand{\thismonth}{\ifcase\month\or
  January\or February\or March\or April\or May\or June\or
  July\or August\or September\or October\or November\or December\fi
  \space\number\year}
\newcommand{\bauth}[1]{\mbox{#1},}
\newcommand{\bart}[1]{\textit{#1},}
\newcommand{\bjourn}[3]{#1 \textbf{#2} (#3)}
\newcommand{\bbook}[1]{\textsl{#1},}

\newcommand{\bpp}[2]{pp.~#1--#2.}
\DeclareMathAlphabet{\mathrmsl}{OT1}{cmr}{m}{sl}
\newcommand{\symb}[2]{\newcommand{#1}{\mathit{#2}} }
\newcommand{\rssymb}[2]{\newcommand{#1}{\mathrmsl{#2}} }
\newcommand{\oper}[3][n]{\newcommand{#2}{\mathop{\mathrm{#3}}%
\ifx n#1\nolimits\else\limits\fi} }
\newcommand{\rsoper}[3][n]{\newcommand{#2}{\mathop{\mathrmsl{#3}}%
\ifx n#1\nolimits\else\limits\fi} }
\oper\End{End}
\oper\Aut{Aut}
\oper\Ad{Ad}
\oper\Vol{Vol}
\oper\Alt{Alt}
\oper\Sym{Sym}
\rsoper{\sym}{sym}
\rsoper{\alt}{alt}
\rsoper\trace{tr}
\rsoper\divg{div}
\symb\vol{vol}
\symb\ev{ev}
\rssymb\iden{id}
\rssymb{\ric}{ric}
\rssymb{\scal}{scal}
\newcommand{\ip}[1]{\langle#1\rangle}
\newcommand{\restr}[1]{|_{\lower.5pt\hbox{${}_{#1}$}}}
\newcommand{\R}{{\mathbb R}}
\newcommand{\C}{{\mathbb C}}

\newcommand{\eps}{\varepsilon}
\renewcommand{\geq}{\geqslant}

\newcommand{\dsum}{\oplus}
\newcommand{\tens}{\mathbin{\otimes}}
\newcommand{\skwend}{\mathinner{\vartriangle}}

\newcommand{\setdif}{\smallsetminus}
\newcommand{\dual}{^{*\!}}
\newcommand{\Cinf}{\mathrm{C}^\infty}

\newcommand{\SU}{\mathrm{SU}}
\newcommand{\connect}{\#}

\newcommand{\cip}{\ip}

\newcommand{\cH}{{\mathcal H}}
\newcommand{\cO}{{\mathcal O}}
\newcommand{\cL}{{\mathcal L}}
\newcommand{\cS}{{\mathcal S}}

\newcommand{\CP}[1]{\C P^{#1}}

\newcommand{\hv}{\mathcal X}
\def\s/{^{\vphantom{,}}_}
\begin{document}
\title{The geometry of the Toda equation}
\author{David M. J. Calderbank}
\address{Department of Mathematics and Statistics\\ University of Edinburgh\\
King's Buildings, Mayfield Road\\ Edinburgh EH9 3JZ\\ Scotland.}
\email{davidmjc@maths.ed.ac.uk}
\date{\thismonth}
\begin{abstract}
I show that solutions of the $\SU(\infty)$ Toda field equation generating a
fixed Einstein-Weyl space are governed by a linear equation on the
Einstein-Weyl space. From this, obstructions to the existence of Toda
solutions generating a given Einstein-Weyl space are found. I also give a
classification of Einstein-Weyl spaces arising from the Toda equation in more
than one way. This classification coincides with a class of spaces found by
Ward and hence clarifies some of their properties. I end by discussing the
simplest examples.
\end{abstract}
\maketitle
\section{Introduction}

In~\cite{Ward}, Ward showed that solutions $u(x,y,z)$ of the $\SU(\infty)$
Toda field equation $u_{xx}+u_{yy}+(e^u)_{zz}=0$ may be used to define three
dimensional Einstein-Weyl spaces. A Weyl space is a conformal manifold $M$
together with a compatible torsion-free connection (called a \emphdef{Weyl
connection}) and it is said to be \emphdef{Einstein-Weyl} iff the symmetric
tracefree part of the Ricci tensor of this connection
vanishes~\cite{Hitchin3}. Weyl connections on a conformal manifold correspond
bijectively to covariant derivatives (called Weyl derivatives) on the density
line bundle $L^1$, which is the oriented real line bundle whose $n$th power is
$|\Lambda^n TM|$ where $n=\dim M$.

A Weyl space may be described by a choice of compatible Riemannian metric $g$
and the connection $1$-form $\omega$ of the Weyl derivative on $L^1$ relative
to the trivialisation of $L^1$ determined by the volume form of the metric (so
that, for the induced Weyl connection, $Dg=-2\omega\tens g$). In these terms,
the Einstein-Weyl space defined by the solution $u$ of the Toda equation may
be written:
\begin{equation}\label{TodaForm}\begin{split}
g &= e^u(dx^2+dy^2)+dz^2\\ \omega &= -u_z\,dz.
\end{split}\end{equation}

The Toda equation is a nonlinear integrable system, but very few solutions are
known explicitly~\cite{BF,CT,GD,Tod2}. Ward found an implicit procedure for
generating a family of solutions from axially symmetric harmonic functions
$V$. The Einstein-Weyl spaces determined by these implicit solutions are
nevertheless completely explicit (in terms of $V$) and Ward suggested that
``\ldots further investigation is needed to clarify the nature and properties
of this family of spaces''~\cite{Ward}.

In this paper I show that these Einstein-Weyl spaces are precisely the
Einstein-Weyl spaces which can be written in the form~\eqref{TodaForm} in at
least two inequivalent ways. The key observation is that the solutions of the
Toda equation on a fixed Einstein-Weyl background are essentially given by
solutions of a linear system in this Weyl geometry.  More precisely, the
solutions of this linear system, the ``Toda structures'', correspond to
solutions of the Toda equation on the Einstein-Weyl space up to changes of
isothermal coordinates $(x,y)$ and translation of $z$. As a consequence, I
show that an Einstein-Weyl space admits at most a four dimensional space of
compatible Toda structures, with equality iff the Einstein-Weyl space is
Einstein. Furthermore, obstructions to the existence of Toda structures on a
given Einstein-Weyl space are found. These obstructions are sufficient to
establish which of the local forms of compact Einstein-Weyl spaces (found
in~\cite{Tod1}) admit Toda structures.

In section~\ref{dftf}, I prove that the existence of more than one Toda
structure on an Einstein-Weyl space is equivalent to the existence of a
conformal vector field of a special type, which I will call an ``axial
symmetry''. This fact is used in section~\ref{ward} to classify the resulting
spaces. The simplest examples are then discussed in the final section.

I work throughout with the density bundles $L^w$ ($w\in\R$). A conformal
structure may then be defined as an $L^2$ valued metric, so that the conformal
inner product of vector fields $X,Y$ is $\cip{X,Y}\in\Cinf(M,L^2)$. Compatible
Riemannian metrics correspond to trivialisations of $L^1$, and such a
trivialisation is often called a \emphdef{length scale} or
\emphdef{gauge}. When tensoring with a density line bundle, I shall omit the
tensor product sign, and sections of $L^{w-1}TM$ or $L^{w+1}T\dual M$ are
called vector fields or $1$-forms of weight $w$ respectively. The Hodge star
operator on an oriented conformal $3$-manifold identifies $L^w$ with
$L^{w+3}\Lambda^3T\dual M$ and $L^{w+1}T\dual M$ with $L^{w+2}\Lambda^2T\dual
M$ and it will be taken to have square $-\iden$. For further details
see~\cite{CP2,Gauduchon4}. The results in this paper are local in character,
and so, where necessary, vector fields are taken to be nonvanishing and
manifolds simply connected.

\subsection*{Acknowledgements} Thanks to Henrik Pedersen and Paul Tod
for interesting and helpful discussions.

\section{Toda structures on Einstein-Weyl spaces}

The Einstein-Weyl spaces arising from the Toda Ansatz~\eqref{TodaForm}
have been characterised by Tod~\cite{Tod2} as those which admit a
shear-free twist-free geodesic congruence. If
$\chi\in\Cinf(M,L^{-1}TM)$ denotes the weightless unit vector
field tangent to this congruence (an oriented foliation with one dimensional
leaves) then this means that
\begin{equation}\label{Chi}
D\chi=\tau(\iden-\cip{\chi,.}\tens\chi)
\end{equation}
where $D$ is the Weyl connection, $\tau$ is a section of $L^{-1}$ and
$\cip{\chi,.}$ denotes the weightless $1$-form dual to $\chi$ with respect to
the conformal structure. In~\eqref{TodaForm}, the congruence generated by
$\partial/\partial z$ has this property and one finds that
$2\tau\cip{\chi,.}=-u_z\,dz$. Hence the Weyl derivative $D-2\tau\cip{\chi,.}$
is induced by the Levi-Civita connection of $g$ and so the metric $g$ is
canonically determined, up to a constant multiple, by the Weyl structure and
the congruence~\cite{CP2,Tod2}. I will denote by $\mu$ the trivialisation of
$L^1$ corresponding to $g$ and refer to this gauge $\mu$ (unique up to a
constant) as the \emphdef{LeBrun-Ward} gauge~\cite{LeBrun1,Ward}. The
congruence $\chi$ determines (in principle) the solution of the Toda equation
up to the
choice of isothermal coordinates $(x,y)$ and affine changes of $z$. Fixing the
LeBrun-Ward gauge $\mu$ determines $z$ up to translation. In~\eqref{TodaForm},
$\tau=-\frac12 u_z\mu^{-1}$, $\chi=\mu^{-1}\partial/\partial z$ and
$\cip{\chi,.}=\mu\,dz$.

The equation~\eqref{Chi} for $\chi$ is apparently nonlinear, but it actually
becomes linear as an equation for the weight $\frac12$ vector field
$\hv=\mu^{1/2}\chi$: explicitly, $D\hv=\sigma\,\iden$, for some section
$\sigma=\mu^{1/2}\tau$ of $L^{-1/2}$. Conversely, if $D\hv$ is a multiple of
the identity, then $\chi=\hv/|\hv|$ is a shear-free twist-free geodesic
congruence and $\mu=|\hv|^2$ is the LeBrun-Ward gauge.

Although this observation is trivial, it is the key idea behind the results of
this paper, and so I should explain its origins. In~\cite{LeBrun1}, LeBrun
gave a characterisation of the Toda Einstein-Weyl spaces in terms of
\emphdef{minitwistor theory}~\cite{Hitchin3}. The space of oriented geodesics
in a three dimensional Einstein-Weyl space is a complex surface $\cS$
containing rational curves (``minitwistor lines'') with normal bundle
$\cO(2)$, and shear-free geodesic congruences correspond to divisors in $\cS$
of degree $2$ on each minitwistor line. LeBrun noticed that if the congruence
is also twist-free, then the corresponding divisor is actually a divisor for
$K_\cS^{-1/2}$, where $K_\cS$ is the canonical bundle of $\cS$. After
incorporating the choice of homothety factor of the LeBrun-Ward gauge, Toda
structures on a fixed Einstein-Weyl space correspond to holomorphic sections
of $K_\cS^{-1/2}$. This immediately suggests that a linear equation is
involved, and by applying the Penrose transform, following Tsai~\cite{Tsai},
one finds that sections of $K_\cS^{-1/2}$ correspond to weight $\frac12$
vector fields with tracelike derivative. It is then not hard to guess the
relationship between such a vector field and $\chi$.

\begin{defn} A \emphdef{Toda structure} on a three dimensional Einstein-Weyl
space is a shear-free twist-free geodesic congruence together with a choice of
homothety factor for the corresponding LeBrun-Ward gauge.
\end{defn}
A Toda structure gives (perhaps only implicitly) a solution of the Toda
equation up to changes of isothermal coordinates $(x,y)$ and translation of
$z$.
\begin{prop} Toda structures correspond to solutions of the following
closed linear system for a nonzero weight $\frac12$ vector field $\hv$ and
a $-\frac12$ density $\sigma$:
\begin{align}\label{eq1}
D\hv&=\sigma\,\iden\\
D\sigma&=-\frac12F^D(\hv,.)-\frac16\scal^D\cip{\hv,.}\label{eq2}
\end{align}
where $D$ is the Weyl connection, $F^D$ is its curvature on $L^1$, and
$\scal^D$ is its scalar curvature, which is a section of $L^{-2}$. Hence Toda
structures are parallel sections with respect to a natural connection on
$L^{-1/2}TM \dsum L^{-1/2}$.
\proofof{prop} Equation~\eqref{eq1} has already been established.
Differentiating it and skew-symmetrising yields
$R^{D,\frac12}_{X,Y}\hv=(D_X\sigma) Y-(D_Y\sigma)X$, where $R^{D,\frac12}$
denotes the curvature of $D$ on $L^{-1/2}TM$. Since $D$ is Einstein-Weyl,
\begin{equation}\label{EWcurv}
R^{D,\frac12}_{X,Y}=-\tfrac16\scal^D\cip{X,.}\skwend Y
+\tfrac12 F^D(X,.)\skwend Y-\tfrac12 F^D(Y,.)\skwend X+\tfrac12 F^D(X,Y)\iden
\end{equation}
where (for any $1$-form $\alpha$ and vector fields $X,Y$)
$\alpha\skwend X(Y)=\alpha(Y)X-\cip{X,Y}\flat\alpha$.
Equation~\eqref{eq2} now follows by taking a trace.
\end{proof}
\begin{cor} An Einstein-Weyl space admits at most a four dimensional
space of Toda structures, and hence at most a three parameter family of
shear-free twist-free geodesic congruences.
\end{cor}
By computing the curvature of the connection
\begin{equation*}
D(\hv,\sigma)=\bigl(D\hv-\sigma\,\iden,D\sigma+\tfrac12F^D(\hv,.)
+\tfrac16\scal^D\cip{\hv,.}\bigr)
\end{equation*}
on $L^{-1/2}TM \dsum L^{-1/2}$, one can find obstructions to the existence of
Toda structures on Einstein-Weyl spaces. In particular, substituting
equation~\eqref{eq2} back into $R^{D,\frac12}_{X,Y}\hv=(D_X\sigma)
Y-(D_Y\sigma)X$, and using~\eqref{EWcurv}, yields
$$F^D(X,Y)\hv+\cip{X,\hv}F^D(Y)-\cip{Y,\hv}F^D(X)=0$$
where $F^D(X)=\flat F^D(X,.)$. This condition on $\hv$ is simply that
$\cip{\hv,.}\wedge F^D=0$, or equivalently, $\cip{\hv,*F^D}=0$.
\begin{prop}\label{orth} The congruence associated to a Toda
structure on an Einstein-Weyl space $(M,D)$ must be orthogonal to $*F^D$.
Hence $M$ admits a four dimensional space of Toda structures if and only if
$F^D=0$, i.e., the Einstein-Weyl space is Einstein.
\end{prop}
\noindent The sufficiency of $F^D=0$ follows by verifying that on each of the
three Einstein spaces, the Toda structures (given by Tod~\cite{Tod2}) do
indeed form a four parameter family (where one of the parameters is
essentially the homothety factor of the Einstein metric).

The remaining curvature obstructions are obtained by differentiating
equation \eqref{eq2} and skew-symmetrising. The resulting constraint
on $(\hv,\sigma)$ is:
\begin{multline*}
(D_XF^D)(Y,\hv)-(D_YF^D)(X,\hv)
+\tfrac13\bigl(\cip{X,\hv}D_Y\scal^D-\cip{Y,\hv}D_X\scal^D\bigr)\\
=F^D(X,Y)\sigma.
\end{multline*}
The Cotton-York curvature of the underlying conformal structure
may be defined by $C_{X,Y}Z=(D_XF^D)(Y,Z)-(D_YF^D)(X,Z)
+\tfrac16\bigl(\cip{X,Z}D_Y\scal^D-\cip{Y,Z}D_X\scal^D\bigr)$.
Hence the above constraint relates $C$ to $D\scal^D$ and $F^D$.
Since it is skew in $X,Y$, it is convenient to apply the star operator
to obtain
$$\mathcal Y(\hv,.)+\tfrac16(*D\scal^D)(\hv,.)=\sigma\,{*F^D}$$ where
$\mathcal Y(U,V)=\cip{*(C_{.,.}U),V}$ (which is well known to define a
symmetric tracefree tensor). One simple consequence of this is the following
refinement of Proposition~\ref{orth}.
\begin{prop}\label{null}
The congruence associated to a Toda structure on an Einstein-Weyl
space $(M,D)$ must be orthogonal to $*F^D$ and null with respect to the
Cotton-York tensor $\mathcal Y$. Hence $M$ can only admit a Toda structure if
$\mathcal Y$ is indefinite on the orthogonal complement of $*F^D$.
\end{prop}
To see that this obstruction is nontrivial, I will apply it in the case that
the Weyl structure is given by $(g,\omega)$ with $\omega$ dual to a Killing
field of $g$. On a compact Einstein-Weyl space, there is a unique compatible
metric (up to a constant) with this property, and the Einstein-Weyl structures
satisfying this condition have been classified~\cite{Tod1}. In order to avoid
a case-by-case computation of $\mathcal Y$, I will derive a general formula.
\begin{prop}\label{formulae}
Suppose $D=D^g+\omega$ is Einstein-Weyl with $\omega$ dual to a
Killing field of $g$. Then
\begin{enumerate}
\item $D^g_XF^D=\frac13\scal^D\,\omega\wedge\cip{X,.}$
and so ${*F^D}$ is also dual to a Killing field of $g$.
\item $\mathcal Y(U,V)=\frac32\bigl(\omega(U)(*F^D)(V)+\omega(V)(*F^D)(U)\bigr)
-\cip{\omega,*F^D}\cip{U,V}$.
\end{enumerate}
\proofof{prop} Since $F^D=d\omega$ is closed and $D^g\omega$ is skew,
$D^g_XF^D(Y,Z)=-2(R^g_{Y,Z}\omega)(X)$. The usual formulae for the Ricci
tensor of $g$~\cite{Gauduchon4,Tod1} yield the first result by direct
calculation.

Next observe that $D\scal^D=D^g\scal^D-2\scal^D\omega$ and that
$D_XF^D(Y,Z)=D^g_XF^D(Y,Z) -F^D(\omega\skwend X (Y),Z)-F^D(Y,\omega\skwend X
(Z))-2\omega(X)F^D(Y,Z)$. Also, by~\cite{Gauduchon4},
$D^g\scal^D=3D^g|\omega|^2$
and $D^g\omega=\frac12F^D$, which leads to the following
formula for $C$:
\begin{align*}
C_{X,Y}=&-\omega(X)F^D(Y,.)+\omega(Y)F^D(X,.)\\
&+\tfrac32F^D(\flat\omega,X)\cip{Y,.}-\tfrac32F^D(\flat\omega,Y)\cip{X,.}
+2F^D(X,Y)\omega.
\end{align*}
Applying the star operator gives the second formula.
\end{proof}
\begin{cor}\label{cobs}
Suppose $D=D^g+\omega$ is Einstein-Weyl on $M$ with $\omega$
dual to a Killing field of $g$. Then $(M,D)$ cannot admit a Toda structure
unless ${*F^D}$ is orthogonal to $\omega$.
\end{cor}
Examining the explicit solutions in~\cite{Tod1}, one can easily determine for
which spaces this holds: in terms of the parameters in Case 1 of~\cite{PT}
(which is the generic case), this condition is $abc=0$. In particular, among
the Berger spheres (given by $b=\pm c$ and $a\neq0$), only the round sphere is
Toda, verifying (in another way) the final remarks of~\cite{Tod2}.

\section{Toda structures and symmetries}\label{dftf}

I turn now to the question: which Einstein-Weyl spaces admit more than a one
dimensional family of Toda structures? In the minitwistor space picture, two
Toda structures correspond to two holomorphic sections of $K_\cS^{-1/2}$. Their
Wronskian, being a section of $K_\cS^{-1}\tens T\dual\cS\cong T\cS$, is a
holomorphic vector field on $\cS$. This symmetry of the minitwistor space
induces a symmetry of the Einstein-Weyl space.

\begin{prop} Suppose $\hv_1$ and $\hv_2$ are the weight $\frac12$ vector fields
of two Toda structures. Then $K=*(\hv_1\wedge\hv_2)$ is a divergence-free
twist-free conformal vector field preserving the Weyl connection.
\proofof{prop} Differentiating $K$ gives
$DK=*(\sigma_1\hv_2-\sigma_2\hv_1)=\frac12d^DK$ where
$D\hv_i=\sigma_i\,\iden$. This is skew and so $K$ is a divergence-free
conformal vector field. Also $\cip{K,.}\wedge d^DK=0$ (since $K$ is orthogonal
to $\hv_1$ and $\hv_2$) so $K$ is twist-free. Finally, to show that $K$
preserves the Weyl connection, it suffices to show that the Lie derivative
$\cL_KD=d\trace DK+F^D(K,.)$ of the Weyl derivative on $L^1$ vanishes.  Now
$*F^D$ is orthogonal to $\hv_1$ and $\hv_2$, so $F^D(K,.)=0$, and $\trace
DK=0$ since $K$ is divergence-free.
\end{proof}

Remarkably, the necessary condition of this proposition is also
sufficient.
\begin{thm} An Einstein-Weyl space has a two dimensional family of Toda
structures if and only if it admits a \textup(nonzero\textup) divergence-free,
twist-free conformal vector field preserving the Weyl connection.
\end{thm}
The necessity is the previous proposition. For the converse, suppose that $K$
is a divergence-free, twist-free conformal vector field preserving the Weyl
connection $D$ of an Einstein-Weyl space. Since the result is local, assume
$K$ is nonvanishing. Then $DK=\alpha\skwend K$, for some $1$-form $\alpha$
with $\alpha(K)=0$. Furthermore, if $D^{|K|}$ is the Weyl derivative
corresponding to the trivialisation of $L^1$ given by the length of $K$, then
$D=D^{|K|}+\alpha$, and so $\cL_K\alpha=0$. Next note that, since $K$ is
twist-free, shear-free and divergence-free, it is surface-orthogonal and the
integral surfaces of $K^\perp$ are totally geodesic. The above theorem is now
an immediate consequence of the following proposition.
\begin{prop} Given $D$, $K$, $\alpha$ as above, the covariant derivative
defined by $D^*_X\hv=D_X\hv-\alpha(\hv)X$ is flat on the bundle of vector
fields of weight $\frac12$ orthogonal to $K$.
\proofof{prop} The curvature of $D^*$ is:
\begin{align*}
R^*_{X,Y}\hv&=\bigl(-\tfrac16\scal^D X\skwend Y+\tfrac12 F^D(X,.)\skwend Y
-\tfrac12 F^D(Y,.)\skwend X+\tfrac12 F^D(X,Y)\iden\bigr)(\hv)\\
&\quad-\bigl((D_X\alpha)(\hv)+\alpha(X)\alpha(\hv)\bigr)Y
+\bigl((D_Y\alpha)(\hv)+\alpha(Y)\alpha(\hv)\bigr)X.
\end{align*}
Now since $K$ is a conformal vector field preserving $D$,
$D_X(DK)=R^D_{X,K}$. Also $DK=\alpha\skwend K$, so
$D_X(DK)=(D_X\alpha+\alpha(X)\alpha)\skwend K$. Contracting with $K$
and using the fact that $\alpha(K)=0$ and $\cL_K\alpha=0$ (i.e.,
$(D_X\alpha)(K)=-\alpha(D_XK)=\cip{K,X}|\alpha|^2$) gives
$$D_X\alpha+\alpha(X)\alpha=\frac12 F^D(X)-\frac16\scal^D X^\perp
+|\alpha|^2X^{||},$$
where $X^{||}$ and $X^\perp$ denote the components of $X$ parallel
and orthogonal to $K$. Substituting this into the formula for $R^*$
gives, for $\hv$ orthogonal to $K$,
$$R^*_{X,Y}\hv=\tfrac12F^D(X,Y)\hv-\tfrac12\cip{\hv,Y}\flat F^D(X,.)
+\tfrac12\cip{\hv,X}\flat F^D(Y,.).$$
This vanishes for all $X,Y$ because $F^D(K,.)=0$ and $\hv$ is orthogonal
to $K$, so $\cip{\hv,.}\wedge F^D=0$.
\end{proof}
The parallel sections of $D^*$ satisfy $D\hv=\alpha(\hv)\iden$ and hence give
a two dimensional family of Toda structures.

A consequence of this theorem is the following converse to
Corollary~\ref{cobs}.
\begin{prop}\label{cistoda}
Suppose $D=D^g+\omega$ is Einstein-Weyl on $M$ with $\omega$
dual to a Killing field of $g$ and that ${*F^D}$ is orthogonal to $\omega$.
Then ${*F^D}$ is dual to a divergence-free twist-free conformal vector
field preserving the Weyl connection, and so $M$ admits a two dimensional
family of Toda structures.
\proofof{prop}
Let $K=\mu_g^3\cip{*F^D,.}$ be the vector field dual to ${*F^D}$ with
respect to $g$, where $\mu_g$ is the trivialisation of $L^1$ determined by
$g$. Now by Proposition~\ref{formulae},
$D^gK=-\frac13\scal^D\mu_g^3\,{*\omega}$
(here $*\omega$ is viewed as a skew endomorphism).
Since $D=D^g+\omega$, $DK=-\frac13\scal^D\mu_g^3\,{*\omega}+\omega\skwend K
+\omega(K)\iden$. Now if ${*F^D}$ is orthogonal to $\omega$ then
$\omega(K)=0$ and $-\frac13\scal^D\mu_g^3{*\omega}=\alpha^g\skwend K$
for some $1$-form $\alpha^g$. Hence $K$ is a divergence-free
twist-free conformal vector field, and it preserves the Weyl
connection since $F^D(K,.)=0$, by definition.
\end{proof}

This result could also be easily established by considering each case
in turn (most of which are straightforward). These spaces will feature
in the final section.

\section{Einstein-Weyl spaces with an axial symmetry}\label{ward}

In this section, I will find explicitly all the Einstein-Weyl spaces admitting
a two dimensional family of Toda structures. According to the previous
section, this is equivalent to classifying the Einstein-Weyl spaces admitting
a divergence-free twist-free conformal vector field $K$ preserving the Weyl
connection. I will say that these spaces are Einstein-Weyl \emphdef{with an
axial symmetry}. On such a space, there is a two dimensional family of Toda
structures given by the weight $\frac12$ vector fields $\hv$ orthogonal to $K$
and satisfying $D\hv=\alpha(\hv)\iden$, where $DK=\alpha\skwend K$. In
particular, $D_K\hv=\alpha(\hv)K=D_\hv K$, so $\cL_K\hv=0$ and these Toda
structures are $K$-invariant.

Pick one such Toda structure $\hv$. Then $\alpha(\hv/|\hv|)$ is the section
$\tau$ of $L^{-1}$ given by this Toda structure, and it is only identically
zero if $\hv$ is a parallel vector field (which can only happen on flat
space).  As shown by LeBrun~\cite{LeBrun1}, $\tau$ is a solution of the
abelian monopole equation and applying the Jones and Tod
construction~\cite{JT} to this solution gives a hyperK\"ahler metric with a
Killing field $X$~\cite{BF,GD}. The Toda structure is $K$-invariant, so $K$
lifts to give an additional Killing field of the hyperK\"ahler metric. Since
$K$ and $X$ commute, some linear combination must be a triholomorphic Killing
field and hence the hyperK\"ahler metric arises via the Gibbons-Hawking
Ansatz~\cite{GH} from a harmonic function on $\R^3$. This harmonic function is
invariant under a Killing field of $\R^3$ and, since $K$ is twist-free, one
readily finds that this Killing field must also be twist-free
(see~\cite{CP2}). Hence it is a rotational vector field, and the harmonic
function is axially symmetric. This proves the following result.
\begin{thm} Let $M$ be Einstein-Weyl with an axial symmetry. Then if
$M$ is not flat \textup(with translational symmetry\textup), it is one of
Ward's Einstein-Weyl spaces constructed from an axially symmetric harmonic
function on $\R^3$~\cite{Ward}, and is therefore given explicitly by:
\begin{align*}
g&=(V_\rho^2+V_\eta^2)(d\rho^2+d\eta^2)+d\psi^2\\
\omega&=\frac{2V_\rho V_\eta\,d\eta+(V_\rho^2-V_\eta^2)d\rho}
{\rho(V_\rho^2+V_\eta^2)}
\end{align*}
where $(\rho V_\rho)_\rho+\rho V_{\eta\eta}=0$.
\end{thm}
Note that the monopole on $\R^3$ is $V_\eta$: the choice of the integral $V$
of $V_\eta$ corresponds to the choice of the quotient of the Gibbons-Hawking
metric (see~\cite{CP2}). This freedom involves adding multiples of $\log\rho$
to $V$. Note also that if $V=\log\rho$, then the monopole $V_\eta$
degenerates, the Einstein-Weyl space above is $\R^3$ itself, and
$\partial/\partial\psi$ is the axial symmetry.

The equation for $V$ may be viewed as an equation on $\cH^2$, by thinking of
$V$ as being in the kernel of the conformal Laplacian on $\R^3\setdif\R$, which
is conformal to $S^1\times\cH^2$. More explicitly, if $v=\rho^{1/2}V$, then
$v_{\rho\rho}+v_{\eta\eta}=-\frac14\rho^{-2}v$ and so $v$ is an eigenfunction
of the Laplacian with eigenvalue $\frac18\scal_{\cH^2}$ on the hyperbolic
$2$-space $\cH^2$ with metric $(d\rho^2+d\eta^2)/\rho^2$.

The original choice of Toda structure $\hv$ may be found by rescaling $g$ by
$\rho^2$ to obtain the Weyl structure in the LeBrun-Ward gauge (again,
see~\cite{CP2}):
\begin{align*}
g\s/{LW}&=\rho^2(V_\rho^2+V_\eta^2)(d\rho^2+d\eta^2)+\rho^2\,d\psi^2\\
&=\rho^2(dV^2+d\psi^2)+(\rho V_\eta\,d\rho-\rho V_\rho\,d\eta)^2\\
\omega\s/{LW}&=-\frac{2V_\eta}{\rho^2(V_\rho^2+V_\eta^2)}
(\rho V_\eta \,d\rho-\rho V_\rho \,d\eta).
\end{align*}
The $1$-form $\rho V_\eta \,d\rho-\rho V_\rho \,d\eta$ is locally exact, and
may be integrated explicitly by writing $V=U_\eta$ with $U$ axially
symmetric and harmonic on $\R^3$, so that $z=-\rho U_\rho$.

The other Toda structures come from the radial congruences on $\R^3$ centred
about points on the axis of symmetry.  A more democratic approach involves the
relationship between these examples and Joyce's construction of torus
symmetric scalar flat K{\"a}hler metrics~\cite{Joyce} from a linear equation
on hyperbolic $2$-space. Indeed this linear equation, given in Proposition
3.2.1 of~\cite{Joyce}, is the Cauchy-Riemann form of the equation for $V$,
tensored trivially with $\R^2$. Ignoring the $\R^2$ tensor factor, simply take
$x_1=\rho$, $x_2=\eta$, $\phi_1=\rho V_\eta$, $\phi_2=-\rho V_\rho$ to see
that Joyce's equation is solved by axially symmetric harmonic
functions. However, the advantage of his approach is that the pair $(\rho
V_\eta,-\rho V_\rho)$ is identified with a section $\Phi$ of a square root of
the canonical bundle of $\cH^2$ satisfying an invariant equation. Now the
Einstein-Weyl structure may be written
\begin{align*}
g&=|\Phi|^2g^{\strut}_{\cH^2}+d\psi^2\\
\omega&=\Phi^2/|\Phi|^2
\end{align*}
and so it does not actually depend upon the choice of coordinates
$(\rho,\eta)$ identifying $\cH^2$ with the upper half plane. Such an
identification is given by a choice of a point at infinity on the hyperbolic
disc and each point in this circle gives a Toda congruence.

Two solutions of Joyce's equation generate a scalar-flat K\"ahler metric with
two Killing fields, and Ward's Einstein-Weyl spaces arise as the quotients by
each of these Killing fields. Joyce finds the solution $V=\log\rho$ (which
generates $\R^3$) and superposes it with its image under isometries of
hyperbolic $2$-space (where these isometries are applied to $\Phi$). In this
way he obtains torus symmetric selfdual conformal structures on $k\CP2$,
generalising (for $k\geq4$) the torus symmetric examples obtained from the
hyperbolic Ansatz of LeBrun~\cite{LeBrun1}.

\section{Examples}

The simplest axially symmetric harmonic functions on $\R^3$ are the constant
functions and the fundamental solutions. The most trivial solution $V_\eta=0$,
$V=\log\rho$ yields $\R^3$. If $V_\eta=b$ or $V_\eta=c/\sqrt{\rho^2+\eta^2}$
then the Gibbons-Hawking metric is $\R^4$ and the triholomorphic Killing field
is an infinitesimal translation or selfdual rotation respectively. Hence the
Einstein-Weyl spaces obtained are the quotients of $\R^4$ by Killing fields
(infinitesimal transrotations or rotations) given in~\cite{PT}.

To obtain more complicated examples, one can take linear combinations of
fundamental solutions and constant solutions. In this way, one can find
the Einstein-Weyl quotients of the Taub-NUT and Eguchi-Hanson metrics, more
or less by direct substitution, although more manageable expressions
are obtained after transforming the $(\rho,\eta)$ coordinates.

The Taub-NUT solutions are given by
$V=a\log\rho+b\eta+c\log\frac{\eta+\sqrt{\rho^2+\eta^2}}\rho$ and it is
convenient to set
$\rho=r\cos\theta, \eta=r\sin\theta$ so that $\rho V_\eta=(br+c)\cos\theta$
and $\rho V_\rho=a-c\sin\theta$. Then
\begin{align*}
g\s/{LW}&=\bigl((br+c)^2\cos^2\theta+(a-c\sin\theta)^2\bigr)
(dr^2+r^2\,d\theta^2)+r^2\cos^2\theta \,d\psi^2\\
\omega\s/{LW}&=-\frac{2(br+c)}
{r\bigl((br+c)^2\cos^2\theta+(a-c\sin\theta)^2\bigr)}
\,d\bigl(-ar\sin\theta+\tfrac12br^2\cos^2\theta+cr\bigr).
\end{align*}
Note that $bc=0$ gives the quotients of $\R^4$ mentioned briefly above.

The Eguchi-Hanson solutions are obtained from $$V=a\log\rho+
\tfrac12(b+c/\eps)\log\tfrac{\eta-\eps+\sqrt{\rho^2+(\eta-\eps)^2}}\rho+
\tfrac12(b-c/\eps)\log\tfrac{\eta+\eps+\sqrt{\rho^2+(\eta+\eps)^2}}\rho,$$
where $\eps^2=\pm1$ (without loss of generality).
When $\eps^2=-1$ this is the
potential for an axially symmetric circle of charge, while $\eps^2=+1$
corresponds to two point sources on the axis of symmetry. These cases are
sometimes referred to as Eguchi-Hanson I and II respectively. The former
is always incomplete, but its Einstein-Weyl quotients are perhaps more
interesting than those of Eguchi-Hanson II.

Coordinates adapted to these geometries are obtained via
$\rho=\sqrt{R^2-\eps^2}\sin\theta$ and $\eta=R\cos\theta$ so that
\begin{align*}
\rho V_\eta&=\frac{(bR+c\cos\theta)\sqrt{R^2-\eps^2}\sin\theta}
{R^2-\eps^2\cos^2\theta}  \\ \tag*{and}
\rho V_\rho&=\frac{a(R^2-\eps^2\cos^2\theta)-b(R^2-\eps^2)\cos\theta
+cR\sin^2\theta}{R^2-\eps^2\cos^2\theta}.
\end{align*}
The Toda structure is now given by:
\begin{align*}
g\s/{LW}&=\bigl((a\cos\theta-b)^2(R^2-\eps^2)+(aR+c)^2\sin^2\theta\bigr)
\left(\frac{dR^2}{R^2-\eps^2}+d\theta^2\right)\\
&\qquad+(R^2-\eps^2)\sin^2\theta \,d\psi^2\\
\omega\s/{LW}&=-\frac{2(bR+c\cos\theta)}
{(a\cos\theta-b)^2(R^2-\eps^2)+(aR+c)^2\sin^2\theta}
\,d\bigl(-aR\cos\theta+bR-c\cos\theta\bigr).
\end{align*}
The family given by $a=0,\eps^2=-1$ also arises as a quotient of the
scalar flat K\"ahler metric on $S^2\times\cH^2$. If we write this
as:
$$g=\frac{dR^2}{R^2+1}+(R^2+1)ds^2+d\theta^2+\sin^2\theta \,d\phi^2$$
then $K=b\partial/\partial s+c\partial/\partial\phi$ is a Killing
field. Coordinates adapted to $K$ are given by $\chi=bs+c\phi,
\psi=b\phi-cs$ so that $K$ is a multiple of $\partial/\partial\chi$ and
the quotient metric $g-g(K,.)/g(K,K)$ is
$$\frac{dR^2}{R^2+1}+d\theta^2+\frac{(R^2+1)\sin^2\theta}
{b^2(R^2+1)+c^2\sin^2\theta}\,d\psi^2.$$
This is the same conformal structure as above, and one readily
checks that the Weyl structures also agree. Now $S^2\times\cH^2$
is conformal to $\R^4\setdif\R$ and so these Weyl structures are globally
defined on $S^3$ for $b\neq0$ (since $\partial/\partial s$ is a
dilation). Hence, as remarked in~\cite{CP2}, these quotients
of $\R^4$ by dilation plus planar rotation are Toda (although the
congruences are not globally defined on $S^3$). Additionally, the calculations
of this section verify explicitly that they are Einstein-Weyl
with an axial symmetry, in accordance with Proposition~\ref{cistoda},
and arise from the Eguchi-Hanson I metrics.

\end{document}